\documentstyle[12pt]{article}
\begin{document}
\newtheorem{proposition}{Proposition}[section]
\renewcommand{\theproposition}{{\thesection.}{\arabic{proposition}}}
\renewcommand{\theequation}{{\thesection}.{\arabic{equation}}}
\def\ra{\rightarrow}
\def\cR{{\cal R}}
\def\la{\leftarrow}
\def\sg{\sigma}
\newcommand{\bnpr}{\begin{proposition}}
\newcommand{\edpr}{\end{proposition}}
\def\th{\theta}
\def\sms{\small}
\def\tl{\tilde}
\def\ls{\large}
\def\ld{\ldots}
\newcommand{\id}{{\rm id}}
\def\vth{\vartheta}
\def\lb{\label}
\newcommand{\bneqn}{\begin{eqnarray}}
\newcommand{\edeqn}{\end{eqnarray}}
\def\NN{{\rm N}\!\!\!\!\!{\rm N}}
\def\ZZ{Z\!\!\!\!\!\!Z}
\def\Zp{Z\!\!\!\!\!\!Z_{\,+}}
\def\Z2{Z\!\!\!\!\!\!Z_{\,2}}
\def\Up{\Upsilon}
\def\Dp{\Delta_{+}}
\newcommand{\edtb}{\end{tabular}}
   \def\Up{\Upsilon}
\def\uD{\underline{\Delta}}
\def\ul{\underline}
\def\ol{\overline}
  \def\ot{\otimes}
\def\uDp{\underline{\Delta}_{+}}
\def\nin{\noindent}
\newcommand{\bntb}{\begin{tabular}}
\title{UNIFIED DESCRIPTION OF QUANTUM AFFINE\\
(SUPER)ALGEBRAS $U_q(A_{1}^{(1)})$\\ 
AND $U_q(C(2)^{(2)})$\footnote{Talk given
by V.N. Tolstoy} }

\author{S.M. KHOROSHKIN$^{1}$, J. LUKIERSKI$^{2}$,
V.N. TOLSTOY$^{3}$}

\maketitle 
\begin{center} 
{$^{1}$Institute of Theoretical and and Experimental Physics \\
117259 Moscow \& Russia (e-mail: khor@heron.itep.ru)}\\

\vskip 0.2cm 
{$^{2}$Institute of Theoretical Physics, University of Wroc\l aw\\
50-204 Wroc\l aw \& Poland (e-mail: lukier@ift.uni.wroc.pl)}\\

\vskip 0.2cm 
{$^3$ Institute of Nuclear Physics, Moscow State University \\
119899 Moscow \& Russia (e-mail: tolstoy@nucl-th.npi.msu.su)}
\end{center} 
\vskip 0.5cm 
\begin{abstract}  
We show that the quantum affine algebra $U_{q}(A_{1}^{(1)})$ and
the quantum affine superalgebra $U_{q}(C(2)^{(2)})$ admit  unified
description. The difference between them 
 consists in the
 phase factor which is equal to $1$ for $U_{q}(A_{1}^{(1)})$
 and  is  equalto $-1$ for $U_{q}(C(2)^{(2)})$. We present such a description
for the construction of Cartan-Weyl generators and their commutation
relations, as well for the universal R-matrices.
\end{abstract}  
\setcounter{equation}{0}
\section{Introduction}
Among variety of all affine Lie (super)algebras\footnote{We introduce
the prefix "super" in brackets to stress that the Lie (super)algebras
include  the Lie algebras as well as the Lie superalgebras.} (both
quantized and non-quantized) the affine (super)algebras of rank 2
play a key role. In the first place, all affine series of the
type $A(n|m)^{(1)}$, $B(n|m)^{(1)}$, $C(n)^{(1)}$, $D(n|m)^{(1)}$,
$A(2n|2m\!-\!1)^{(2)}$, $A(2n\!-\!1|2m\!-\!1)^{(2)}$,$C(n)^{(2)}$,
$D(n|m)^{(2)}$ and $A(2n|2m)^{(4)}$ begin with the affine
(super)algebras of rank 2. Secondly, the contragredient Lie
(super)algebras of rank 2 are basic structural blocks of any
affine (super)algebras of arbitrary rank. This fact permits,
for example,
 the reduction of the proofs of basic theorems for the extremal
projector and the universal R-matrix to the proofs of such theorems
for the (super)algebras of rank 2 (see Refs.~\cite{AST}, 
\cite{T1} - \cite{TK}, \cite{KT1} - \cite{KT4}).
  Further, the
representation theory of the affine (super)algebras (both quantized
and non-quantized) contains some typical elements of the representation
theory of the affine (super)algebras of rank 2.
Besides, in applications of the affine
(super)algebras, first of all the affine (super)algebras of rank 2
are used by virtue of their simplicity.

In this paper along the line of considerations presented in
\cite{nieznany} for rank 2 affine superalgebra
 $U_{q}(B(0,1)^{(1)})$ we
 give detailed description of the quantum untwisted affine
algebra $U_{q}(A_{1}^{(1)})$ ($\simeq U_{q}(\widehat{sl}(2))$)
and the quantum twisted affine superalgebra
$U_{q}(C(2)^{(2)})$. 
Moreover our goal is to show that these
quantum (super)algebras are described in unified way.
Namely, we present in unified way their defining relations, 
 the construction of the Cartan-Weyl bases,
the complete list of all commutation relations of the
Cartan-Weyl generators corresponding to all root vectors
 and finally the unified form of their universal
R-matrices. Difference between both considered
quantum (super)algebras is only determined  by a phase factor which is
equal to $1$ for $U_{q}(A_{1}^{(1)})$ and it is equal to $-1$ for
$U_{q}(C(2)^{(2)})$.
This situation is similar to the finite-dimensional case. Namely, in
the paper \cite{KT1} it was shown that all quantum (super)algebras
$U_q(g)$, where $g$ are the finite-dimensional contragredient Lie
(super)algebras of rank 2, are divided into three classes. Each such
class is characterized by the same Dynkin diagram and has the same reduced
root system, provided that we neglect the type (colour) of the roots
(white, grey or black). Consequently, all the
(super)algebras of the same class have  unified defining relations,
unified construction of the Cartan-Weyl basis and its properties,
 as well as unified universal R-matrix. Difference between the (super)algebras
of the same class is determined by some phase factor which takes values
$\pm1$ depending on the colour of the nodes of the Dynkin diagram.

Basic information about the (super)algebras $A_{1}^{(1)}$ and
$C(2)^{(2)}$ is presented in the tables 1a and 1b
(see Refs.~\cite{K1,K2,VdL}). In the table 1a there are listed 
the standard and
symmetric Cartan matrices $A$ and $A^{sym}$, the corresponding
extended symmetric matrices $\bar{A}^{sym}$ and their inverses
$(\bar{A}^{sym})^{-1}$, and also the sets of odd roots (odd roots),
the Dynkin diagrams (diagram), and the dimensions of these
(super)algebras (dim).
We remind some elementary definitions of the colour of the roots:
\begin{itemize}
\item All even roots are called white roots. A white root is
pictured by the white node
\mbox{\begin{picture}(10,10)\put(4,3){\circle{8}}\end{picture}}.
\item An odd root $\gamma$ is called a grey root if
$2\gamma$ is not a root. Such an odd root is pictured by the grey
node $\ot$.
\item An odd root $\gamma$ is called a dark root if
$2\gamma$ is a root. Such an odd root is pictured by the dark node
\mbox{\begin{picture}(10,10)\put(4,3){\circle*{8}}\end{picture}}.
\end{itemize}
We also remind the definition of the reduced system of the positive
root system $\Delta_+$ for any contragredient (super)algebras of finite
growth.
\begin{itemize}
\item The system $\ul{\Delta}_+$ is called the reduced system
if it is defined by the following way:
$\ul{\Delta}_+\!=\Delta_+\!\backslash\{2\gamma\in \Delta_+|\gamma
\;{\rm is\; odd}\}$.
That is the reduced system $\ul{\Delta}_+$ is obtained from the total
system $\Delta_+$ by removing of all doubled roots $2\gamma$ where $\gamma$ is
a dark odd root.
\end{itemize}
The total and reduced root systems of the (super)algebras $A_{1}^{(1)}$
and $C(2)^{(2)}$ are listed in the table 1b. It is convenient to
present the total $\Delta\!=\!\Delta_+\!\bigcup(-\Delta_+)$ and reduced
$\uD\!=\!\uDp\!\bigcup(-\uDp)$ root systems by the pictures: Figs. 1, 2a, 2b.
Comparing Fig. 1 and Fig. 2b we see that the reduced root systems of
$A_{1}^{(1)}$ and $C(2)^{(2)}$ coincide if we neglect colour of the
roots.
\newpage
\vskip 10pt
\centerline{\bf Table 1a}
\vskip 10pt
\nin
{\footnotesize
{\renewcommand{\arraystretch}{0}
\renewcommand{\tabcolsep}{4pt}
\bntb{|lccccc|}
\hline
\rule{0pt}{8pt}&&&&&\\
\strut
{ $g(A,\Up)$}
 & $A=A^{sym}$ & $\bar{A}^{sym}$
&$(\bar{A}^{sym})^{-1}$ & odd & diagram \\[20pt]
\hline
\rule{0pt}{8pt}&&&&&\\
\strut
{ $A_{1}^{(1)}$}
 & $\left(\hspace{-3mm}
\begin{array}{rr}\strut 2&-2\\
\strut -2&2 \end{array}\right)$
& $\left(\begin{array}{rrr}\strut 0&1&0\\
\strut 1&2&-2\\
\strut 0&-2&2\end{array}\right)$ &
$\left(\begin{array}{rrr}\strut 0&1&1\\
\strut 1&0&0\\
\strut 1&0&\frac{1}{2}\end{array}\right)$ & $\emptyset$
& \makebox{\begin{picture}(10,10)\thicklines
\put(-10,3){\circle{8}}
\put(-20,12){\footnotesize$\delta\!-\!\alpha$}
\put(-7,5){\line(1,0){19}}
\put(-7,1){\line(1,0){19}}
\put(16,3){\circle{8}}
\put(14,12){\footnotesize $\alpha$}\end{picture}}\\[25pt]
\strut
{  $C(2)^{(2)}$}
 & $\left(\hspace{-3mm}
\begin{array}{rr}\strut 2&-2\\
\strut -2&2 \end{array}\right)$
& $\left(\begin{array}{rrr}\strut 0&1&0\\
\strut 1&2&-2\\
\strut 0&-2&2\end{array}\right)$ &
$\left(\begin{array}{rrr}\strut 0&1&1\\
\strut 1&0&0\\
\strut 1&0&\frac{1}{2}\end{array}\right)$ & $\{\delta\!-\!\alpha,a\}$
& \makebox{\begin{picture}(10,10)\thicklines
\put(-10,3){\circle*{8}}
\put(-20,12){\footnotesize$\delta\!-\!\alpha$}
\put(-7,5){\line(1,0){19}}
\put(-7,1){\line(1,0){19}}
\put(16,3){\circle*{8}}
\put(14,12){\footnotesize $\alpha$}\end{picture}}\\[25pt]
\hline
\edtb}}
\vskip 20pt
\centerline{\bf Table 1b}
\vskip 10pt
\nin
{\footnotesize
{\renewcommand{\arraystretch}{0}
\renewcommand{\tabcolsep}{4pt}
\bntb{|lcc|}
\hline
\rule{0pt}{8pt}&&\\
\strut $g(A,\Up)$\hspace{1.1cm} & $\Dp$ \hspace{1cm} & $\uDp$
\\[14pt]
\hline
\rule{0pt}{8pt}&&\\
\strut $A_{1}^{(1)}$ \hspace{1.1cm}&
$\{\alpha,\,n\delta\!\pm\!\alpha,\,n\delta\,|\,n\in\NN\}$
\hspace{1cm} & $\{\alpha,\,n\delta\!\pm\!\alpha,\,n\delta\,|\,n\in\NN\}$
\\[20pt]
\strut $C(2)^{(2)}$ \hspace{1.1cm}&
$\{\alpha,\,2\alpha,\,n\delta\!\pm\!\alpha,
\,2n\delta\!\pm\!2\alpha,\,n\delta\,|\,n\in\NN\}$
\hspace{1cm} & $\{\alpha,\,n\delta\!\pm\!\alpha,\,n\delta\,|\,n\in\NN\}$
\\[20pt]
\hline
\edtb}}
\vskip 20pt
\leftline{
\setlength\unitlength{0.26mm}
\begin{picture}(400,100)
\put(0,72){\ldots}
\put(24,78){\footnotesize$-\!4\delta\!+\!\alpha$}
\put(44,72){\circle{5}}
\put(68,78){\footnotesize$-\!3\delta\!+\!\alpha$}
\put(88,72){\circle{5}}
\put(112,78){\footnotesize$-\!2\delta\!+\!\alpha$}
\put(132,72){\circle{5}}
\put(161,78){\footnotesize$-\!\delta\!+\!\alpha$}
\put(176,72){\circle{5}}
\put(218,78){\footnotesize$\alpha$}
\put(220,72){\circle{5}}
\put(255,78){\footnotesize$\delta\!+\!\alpha$}
\put(264,72){\circle{5}}
\put(294,78){\footnotesize$2\delta\!+\!\alpha$}
\put(308,72){\circle{5}}
\put(338,78){\footnotesize$3\delta\!+\!\alpha$}
\put(352,72){\circle{5}}
\put(382,78){\footnotesize$4\delta\!+\!\alpha$}
\put(396,72){\circle{5}}
\put(420,72){\ldots}
\put(0,50){\ldots}
\put(35,56){\footnotesize$-\!4\delta$}
\put(44,50){\circle{5}}
\put(79,56){\footnotesize$-\!3\delta$}
\put(88,50){\circle{5}}
\put(123,56){\footnotesize$-\!2\delta$}
\put(132,50){\circle{5}}
\put(169,56){\footnotesize$-\!\delta$}
\put(176,50){\circle{5}}
\put(220,50){\vector(0,1){19}}
\put(220,50){\vector(2,-1){42}}
\put(262,56){\footnotesize$\delta$}
\put(264,50){\circle{5}}
\put(303,56){\footnotesize$2\delta$}
\put(308,50){\circle{5}}
\put(347,56){\footnotesize$3\delta$}
\put(352,50){\circle{5}}
\put(391,56){\footnotesize$4\delta$}
\put(396,50){\circle{5}}
\put(4200,50){\ldots}
\put(0,28){\ldots}
\put(24,33){\footnotesize$-\!4\delta\!-\!\alpha$}
\put(44,28){\circle{5}}
\put(68,33){\footnotesize$-\!3\delta\!-\!\alpha$}
\put(88,28){\circle{5}}
\put(112,33){\footnotesize$-\!2\delta\!-\!\alpha$}
\put(132,28){\circle{5}}
\put(161,33){\footnotesize$-\!\delta\!-\!\alpha$}
\put(176,28){\circle{5}}
\put(213,33){\footnotesize$-\!\alpha$}
\put(220,28){\circle{5}}
\put(255,33){\footnotesize$\delta\!-\!\alpha$}
\put(264,28){\circle{5}}
\put(294,33){\footnotesize$2\delta\!-\!\alpha$}
\put(308,28){\circle{5}}
\put(338,33){\footnotesize$3\delta\!-\!\alpha$}
\put(352,28){\circle{5}}
\put(382,33){\footnotesize$4\delta\!-\!\alpha$}
\put(396,28){\circle{5}}
\put(420,28){\ldots}
\end{picture}
}
\vskip -14pt
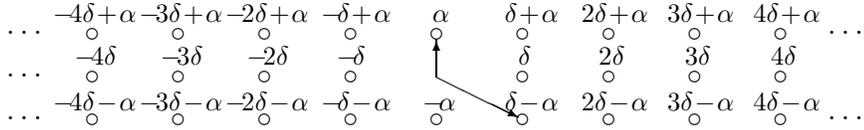
\centerline{\footnotesize Fig. 1. The total and
reduced root system ($\Delta=\uD$) of $A_{1}^{(1)} (\simeq
\widehat{sl}_{2})$.}
\vskip 20pt
\leftline{
\setlength\unitlength{0.26mm}
\begin{picture}(400,120)
\put(0,94){\ldots}
\put(24,101){\footnotesize$-\!4\delta\!+\!2\alpha$}
\put(44,94){\circle{5}}
\put(112,101){\footnotesize$-\!2\delta\!+\!2\alpha$}
\put(132,94){\circle{5}}
\put(215,101){\footnotesize$2\alpha$}
\put(220,94){\circle{5}}
\put(294,101){\footnotesize$2\delta\!+\!2\alpha$}
\put(308,94){\circle{5}}
\put(382,101){\footnotesize$4\delta\!+\!2\alpha$}
\put(396,94){\circle{5}}
\put(420,94){\ldots}
\put(0,72){\ldots}
\put(24,78){\footnotesize$-\!4\delta\!+\!\alpha$}
\put(44,72){\circle*{5}}
\put(68,78){\footnotesize$-\!3\delta\!+\!\alpha$}
\put(88,72){\circle*{5}}
\put(112,78){\footnotesize$-\!2\delta\!+\!\alpha$}
\put(132,72){\circle*{5}}
\put(161,78){\footnotesize$-\!\delta\!+\!\alpha$}
\put(176,72){\circle*{5}}
\put(218,78){\footnotesize$\alpha$}
\put(220,72){\circle*{5}}
\put(255,78){\footnotesize$\delta\!+\!\alpha$}
\put(264,72){\circle*{5}}
\put(294,78){\footnotesize$2\delta\!+\!\alpha$}
\put(308,72){\circle*{5}}
\put(338,78){\footnotesize$3\delta\!+\!\alpha$}
\put(352,72){\circle*{5}}
\put(382,78){\footnotesize$4\delta\!+\!\alpha$}
\put(396,72){\circle*{5}}
\put(420,72){\ldots}
\put(0,50){\ldots}
\put(35,56){\footnotesize$-\!4\delta$}
\put(44,50){\circle{5}}
\put(79,56){\footnotesize$-\!3\delta$}
\put(88,50){\circle{5}}
\put(123,56){\footnotesize$-\!2\delta$}
\put(132,50){\circle{5}}
\put(169,56){\footnotesize$-\!\delta$}
\put(176,50){\circle{5}}
\put(220,50){\vector(0,1){19}}
\put(220,50){\vector(2,-1){42}}
\put(262,56){\footnotesize$\delta$}
\put(264,50){\circle{5}}
\put(303,56){\footnotesize$2\delta$}
\put(308,50){\circle{5}}
\put(347,56){\footnotesize$3\delta$}
\put(352,50){\circle{5}}
\put(391,56){\footnotesize$4\delta$}
\put(396,50){\circle{5}}
\put(4200,50){\ldots}
\put(0,28){\ldots}
\put(24,33){\footnotesize$-\!4\delta\!-\!\alpha$}
\put(44,28){\circle*{5}}
\put(68,33){\footnotesize$-\!3\delta\!-\!\alpha$}
\put(88,28){\circle*{5}}
\put(112,33){\footnotesize$-\!2\delta\!-\!\alpha$}
\put(132,28){\circle*{5}}
\put(161,33){\footnotesize$-\!\delta\!-\!\alpha$}
\put(176,28){\circle*{5}}
\put(213,33){\footnotesize$-\!\alpha$}
\put(220,28){\circle*{5}}
\put(255,33){\footnotesize$\delta\!-\!\alpha$}
\put(264,28){\circle*{5}}
\put(294,33){\footnotesize$2\delta\!-\!\alpha$}
\put(308,28){\circle*{5}}
\put(338,33){\footnotesize$3\delta\!-\!\alpha$}
\put(352,28){\circle*{5}}
\put(382,33){\footnotesize$4\delta\!-\!\alpha$}
\put(396,28){\circle*{5}}
\put(420,28){\ldots}
\put(0,6){\ldots}
\put(24,11){\footnotesize$-\!4\delta\!-\!2\alpha$}
\put(44,6){\circle{5}}
\put(112,11){\footnotesize$-\!2\delta\!-\!2\alpha$}
\put(132,6){\circle{5}}
\put(212,11){\footnotesize$-\!2\alpha$}
\put(220,6){\circle{5}}
\put(294,11){\footnotesize$2\delta\!-\!2\alpha$}
\put(308,6){\circle{5}}
\put(382,11){\footnotesize$4\delta\!-\!2\alpha$}
\put(396,6){\circle{5}}
\put(420,6){\ldots}
\end{picture}
}
\vskip 5pt
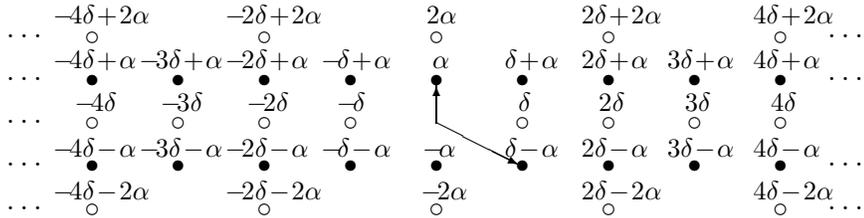
\centerline{\footnotesize Fig. 2a. The total root system $\Delta$
of $C(2)^{(2)}$.}
\vskip 20pt
\leftline{
\setlength\unitlength{0.26mm}
\begin{picture}(400,120)
\put(0,72){\ldots}
\put(24,78){\footnotesize$-\!4\delta\!+\!\alpha$}
\put(44,72){\circle*{5}}
\put(68,78){\footnotesize$-\!3\delta\!+\!\alpha$}
\put(88,72){\circle*{5}}
\put(112,78){\footnotesize$-\!2\delta\!+\!\alpha$}
\put(132,72){\circle*{5}}
\put(161,78){\footnotesize$-\!\delta\!+\!\alpha$}
\put(176,72){\circle*{5}}
\put(218,78){\footnotesize$\alpha$}
\put(220,72){\circle*{5}}
\put(255,78){\footnotesize$\delta\!+\!\alpha$}
\put(264,72){\circle*{5}}
\put(294,78){\footnotesize$2\delta\!+\!\alpha$}
\put(308,72){\circle*{5}}
\put(338,78){\footnotesize$3\delta\!+\!\alpha$}
\put(352,72){\circle*{5}}
\put(382,78){\footnotesize$4\delta\!+\!\alpha$}
\put(396,72){\circle*{5}}
\put(420,72){\ldots}
\put(0,50){\ldots}
\put(35,56){\footnotesize$-\!4\delta$}
\put(44,50){\circle{5}}
\put(79,56){\footnotesize$-\!3\delta$}
\put(88,50){\circle{5}}
\put(123,56){\footnotesize$-\!2\delta$}
\put(132,50){\circle{5}}
\put(169,56){\footnotesize$-\!\delta$}
\put(176,50){\circle{5}}
\put(220,50){\vector(0,1){19}}
\put(220,50){\vector(2,-1){42}}
\put(262,56){\footnotesize$\delta$}
\put(264,50){\circle{5}}
\put(303,56){\footnotesize$2\delta$}
\put(308,50){\circle{5}}
\put(347,56){\footnotesize$3\delta$}
\put(352,50){\circle{5}}
\put(391,56){\footnotesize$4\delta$}
\put(396,50){\circle{5}}
\put(4200,50){\ldots}
\put(0,28){\ldots}
\put(24,33){\footnotesize$-\!4\delta\!-\!\alpha$}
\put(44,28){\circle*{5}}
\put(68,33){\footnotesize$-\!3\delta\!-\!\alpha$}
\put(88,28){\circle*{5}}
\put(112,33){\footnotesize$-\!2\delta\!-\!\alpha$}
\put(132,28){\circle*{5}}
\put(161,33){\footnotesize$-\!\delta\!-\!\alpha$}
\put(176,28){\circle*{5}}
\put(213,33){\footnotesize$-\!\alpha$}
\put(220,28){\circle*{5}}
\put(255,33){\footnotesize$\delta\!-\!\alpha$}
\put(264,28){\circle*{5}}
\put(294,33){\footnotesize$2\delta\!-\!\alpha$}
\put(308,28){\circle*{5}}
\put(338,33){\footnotesize$3\delta\!-\!\alpha$}
\put(352,28){\circle*{5}}
\put(382,33){\footnotesize$4\delta\!-\!\alpha$}
\put(396,28){\circle*{5}}
\put(420,28){\ldots}
\end{picture}
}
\vskip -10pt
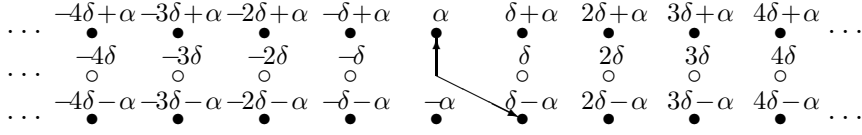
\centerline{\footnotesize Fig. 2b. The reduced root
system $\uD$ of $C(2)^{(2)}$.}

\setcounter{equation}{0}
\section{Defining Relations of $U_{\!q}(\!A_1^{(\!1\!)})$ and
$U_{\!q}(\!C(2)^{(\!2\!)})$}
The quantum (q-deformed) affine (super)algebras $U_q(A_1^{(1)})$
and $U_q(C(2)^{(2)})$ are generated by the Chevalley elements
$k_{\rm d}^{\pm1}:=q^{\pm h_{\rm d}}$,
$k_{\alpha}^{\pm1}:=q^{\pm h_{\alpha}}$, 
$k_{\delta-\alpha}^{\pm1}:=q^{\pm h_{\delta-\alpha}}$,
$e_{\pm\alpha}$, $e_{\pm(\delta-a)}$ with the defining
 relations:  
\bneqn
k_{\gamma}^{}k_{\gamma}^{-1}&\!\!=\!\!&
k_{\gamma}^{-1}k_{\gamma}^{}=1~,\qquad\quad\;
[k_{\gamma}^{\pm 1},k_{\gamma'}^{\pm 1}]=0~,\qquad
\lb{DR1}
\\[7pt]
k_{\gamma}^{}e_{\pm\alpha}^{}k_{\gamma}^{-1}&\!\!=\!\!&
q^{\pm(\gamma,\alpha)}e_{\pm\alpha}^{}~,
\qquad k_{\gamma}^{}e_{\pm(\delta-\alpha)}^{}k_{\gamma}^{-1}=
q^{\pm(\gamma,\delta-\alpha)}e_{\pm(\delta-\alpha)}^{}~,
\lb{DR2}
\edeqn
\bneqn
[e_{\alpha}^{},e_{-\delta+\alpha}^{}]&\!\!=\!\!&0~,\qquad\qquad\qquad\quad
[e_{-\alpha}^{},e_{\delta-\alpha}^{}]=0~,\phantom{possible, poss}
\lb{DR3}
\\[7pt]
[e_{\alpha}^{},e_{-\alpha}^{}]&\!\!=\!\!&[h_{\alpha}]_q~,\qquad\qquad\;\;\,
[e_{\delta-\alpha}^{},e_{-\delta+\alpha}^{}]=[h_{\delta-\alpha}]_q~,
\lb{DR4}
\edeqn
\bneqn
[e_{\pm\alpha}^{},[e_{\pm\alpha}^{},[e_{\pm\alpha}^{},
e_{\pm(\delta-\alpha)}^{}]_q]_q]_{q}&\!\!=\!\!&0~,\quad\;\;
\lb{DR5}
\\[7pt]
[[[e_{\pm\alpha},e_{\pm(\delta-\alpha)}]_q,e_{\pm(\delta-\alpha)}]_q,
e_{\pm(\delta-\alpha)}]_q&\!\!=\!\!&0~,\quad\;\;
\lb{DR6}
\edeqn
where ($\gamma\!=\!{\rm d},\alpha,\delta-\alpha$), $({\rm d},\alpha)\!=\!0$,
$({\rm d},\delta)\!=\!1$, and
$[h_\beta]_q\!:=\!(k_\beta\!-\!k_\beta^{-\!1})/(q\!-\!q^{-\!1})$.
The brackets $[\cdot,\cdot]$ and $[\cdot,\cdot]_{q}$ are the super-,
and q-super-commutators:
\begin{equation}
\begin{array}{rcl}
[e_{\beta}^{},e_{\beta'}^{}]&\!\!=\!\!&e_{\beta}^{}e_{\beta'}^{}-
(-1)^{\vth(\beta)\vth(\beta')}e_{\beta'}^{}e_{\beta}^{}~,
\\[7pt]
[e_{\beta}^{},e_{\beta'}^{}]_{q}&\!\!=\!\!&e_{\beta}^{}e_{\beta'}^{}-
(-1)^{\vth(\beta)\vth(\beta')}q^{(\beta,\beta')}e_{\beta'}^{}e_{\beta}^{}~.
\lb{DR7}
\end{array}
\end{equation}
Here the symbol $\vth(\cdot)$ means the parity function: $\vth(\beta)=0$
for any even root $\beta$ and $\vth(\beta)=1$ for any odd root $\beta$.

{\it Remark}. The left-side sides of the relations (\ref{DR5}) and
(\ref{DR6}) are invariant with respect to the replacement of $q$ by $q^{-1}$.
Indeed, if we remove the q-brackets we see that the left-hand of
(\ref{DR5}) and (\ref{DR6}) contain the symmetric functions of $q$
and $q^{-1}$. This property permits to write the q-commutators in
(\ref{DR5}) and (\ref{DR6}) in the inverse order, i.e.
\bneqn
[[[e_{\pm(\delta-\alpha)}^{},e_{\pm\alpha}^{}]_{q},e_{\pm\alpha}^{}]_{q},
e_{\pm\alpha}^{}]_q&\!\!=\!\!&0\, ,
\lb{DR8}
\\[7pt]
[e_{\pm(\delta-\alpha)},[e_{\pm(\delta-\alpha)},[e_{\pm(\delta-\alpha)},
e_{\pm\alpha}]_q]_q]_q&\!\!=\!\!&0\, .
\lb{DR9}
\edeqn

The standard Hopf structure of the quantum (super)algebras
$U_q(A_1^{(1)})$ and $U_q(C(2)^{(2)})$ is given by the following
formulas for the comultiplication $\Delta_q$ and antipode $S_q$:
\begin{equation}
\begin{array}{rcccl}
\Delta_{q}(k_\gamma^{\pm1})&\!\!=\!\!&k_\gamma^{\pm1}\ot k_\gamma^{\pm1}~,
\qquad\qquad\qquad S_{q}(k_\gamma^{\pm1})&\!\!=\!\!&k_\gamma^{\mp1}~,
\\[7pt]
\Delta_{q}(e_{\beta}^{})&\!\!=\!\!&e_{\beta}^{}\ot 1
+ k_{\beta}^{-1}\ot e_{\beta}^{}~,
\qquad\quad S_{q}(e_{\beta}^{})&\!\!=\!\!&-k_{\beta}e_{\beta}^{}~,
\\[7pt]
\Delta_{q}(e_{-\beta}^{})\!\!&=\!\!&
e_{-\beta}^{}\ot k_{\beta}+1 \ot e_{-\beta}^{}~,
\qquad S_{q}(e_{-\beta}^{})&\!\!=\!\!&-e_{-\beta}^{}k_{\beta}^{-1}~,
\lb{DR15}
\end{array}
\end{equation}
where $\beta=\alpha,\,\delta-\alpha$; $\gamma={\rm d},\,\beta$.
It is not hard to verify by direct calculations for the
defining relations (\ref{DR1})-(\ref{DR6}) that the quantum affine
(super)algebras $U_q(A_1^{(1)})$ and $U_q(C(2)^{(2)})$ have
the following simple involutive (anti)\-auto\-mor\-phisms:

\nin
{\it (i) The non-graded antilinear antiinvolution or conjugation
"$^{*}$"}:
\begin{equation}
\begin{array}{rcccl}
(q^{\pm1})^{*}\!\!&=\!\!&q^{\mp1}~, \qquad\quad
(k_{\gamma}^{\pm 1})^{*}\!\!&=\!\!&k_{\gamma}^{\mp 1}~,\qquad
\\[7pt]
e_{\beta}^{*}\!\!&=\!\!&e_{-\beta}~,\qquad\qquad
e_{-\beta}^{*}\!\!&=\!\!&e_{\beta}
\lb{DR16}
\end{array}
\end{equation}
($(xy)^*=y^*x^*$ for $\forall\;\, x,y\in U_{q}(g)$).

\nin
{\it (ii) The graded antilinear antiinvolution or graded conjugation
"$^{\ddagger}$"}:
\begin{equation}
\begin{array}{rcccl}
(q^{\pm1})^{\ddagger}\!\!&=\!\!&q^{\mp1}~, \qquad\qquad
(k_{\gamma}^{\pm 1})^{\ddagger}\!\!&=\!\!&k_{\gamma}^{\mp 1}~,\qquad
\\[7pt]
e_{\beta}^{\ddagger}\!\!&=\!\!&(-1)^{\vth(\beta)}e_{-\beta}~,\quad\;
e_{-\beta}^{\ddagger}\!\!&=\!\!&e_{\beta}
\lb{DR17}
\end{array}
\end{equation}
($(xy)^{\ddagger}=(-1)^{\deg x\,\deg y} y^{\ddagger}x^{\ddagger}$
for any homogeneous elements $x,y\in U_q(g)$).

\nin
{\it (iii) The Chevalley graded involution $\omega$}:
\begin{equation}
\begin{array}{rcccl}
{}\qquad\omega(q^{\pm1})&\!\!=\!\!&q^{\mp1}~,\qquad\quad
\omega(k_{\gamma}^{\pm1})&\!\!=\!\!&k_{\gamma}^{\pm1}~,
\\[7pt]
\omega(e_{\beta})&\!\!=\!\!&-e_{-\beta}~,\qquad\;
\omega(e_{-\beta})&\!\!=\!\!&-(-1)^{\theta(\beta)}e_{\beta}~.
\lb{DR18}
\end{array}
\end{equation}

\nin
{\it (iv) The Dynkin involution $\tau$} which is associated with the
automorphism of the Dynkin diagrams of the (super)algebras
$A_1^{(1)}$ and $C(2)^{(2)}$:
\begin{equation}
\begin{array}{rcccl}
\tau(q^{\pm1})&\!\!=\!\!&q^{\pm1}~,\qquad\quad
\tau(k_{\rm d}^{\pm1})&\!\!=\!\!&k_{\rm d}^{\pm1}~,
\\[7pt]
\tau(k_{\beta}^{\pm1})&\!\!=\!\!&k_{\delta-\beta}^{\pm1}~,\qquad\;
\tau(k_{-\beta}^{\pm1})&\!\!=\!\!&k_{-\delta+\beta}^{\pm1}~,
\\[7pt]
\tau(e_{\beta})&\!\!=\!\!&e_{\delta-\beta}~,\qquad\;
\tau(e_{-\beta})&\!\!=\!\!&e_{-\delta+\beta}~.
\lb{DR19}
\end{array}
\end{equation}

Here in (\ref{DR16})-(\ref{DR19}) $\beta\!=\!\alpha,\delta\!-\!\alpha$;
$\gamma\!=\!{\rm d},\beta$.

It should be noted that the graded conjugation
"$^{\ddagger}$" and the Chevalley graded involution $\omega$ are
involutive (anti)automorphism of the fourth order, i.e., for example,
$(\omega)^4\!=\!\id$. Note also that the Dynkin involution $\tau$
commutes with all other three involutions, i.e. $\tau(x^*)\!=\!(\tau(x))^*$,
$\tau(x^{\ddagger})\!=\!(\tau(x))^{\ddagger}$ 
and $\omega\tau(x)\!=\!\tau\omega(x)$
for any element $x\in U_q(g)$ ($g\!=\!A_1^{(1)},C(2,0)^{(2)}$).

In the next Section we construct the Cartan-Weyl basis and describe
its properties in detail.

\setcounter{equation}{0}
\section{Cartan-Weyl Basis for $U_{\!q}(\!A_1^{(\!1\!)})$ and
$U_{\!q}(\!C(2)^{(\!2\!)})$}
A general scheme for construction of a Cartan-Weyl basis for
quantized Lie algebras and superalgebras was proposed in Ref.~\cite{T1}
The scheme was applied in detail at first for quantized
finite-dimensional Lie (super)algebras~\cite{KT1} and then to
quantized untwisted affine algebras~\cite{TK}.
This procedure is based on a notion of ``normal ordering'' for the
reduced positive root system. For affine Lie (super)algebras this
notation was introduced in Ref.~\cite{T2}
 (see also Refs.~\cite{T1,KT2,KT3}).

In our case the reduced positive system has only two normal orderings:
\begin{equation}
\mbox{\sms$\alpha,\delta\!+\!\alpha,2\delta\!
+\!\alpha,\ldots,\infty\delta\!+\!\alpha,
\delta,2\delta,3\delta,\ld,\infty\delta,\infty\delta\!-\!\alpha,\ld,
2\delta\!-\!\alpha,\delta\!-\!\alpha$},
\lb{CW1}
\end{equation}
\begin{equation}
\mbox{\sms$\delta\!-\!\alpha,2\delta\!-\!\alpha,\ld,
\infty\delta\!-\!\alpha,\delta,2\delta,3\delta,\ldots,\infty\delta,
\infty\delta\!+\!\alpha,\ld,2\delta\!+\!\alpha,\delta\!+\!\alpha,\alpha$}.
\lb{CW2}
\end{equation}
\nin
The first normal ordering (\ref{CW1}) corresponds to ``clockwise''
ordering for positive roots in Fig. 1, 2b if we start from root $\alpha$
to root $\delta\!-\!\alpha$. The inverse normal ordering (\ref{CW2})
corresponds to ``anticlockwise'' ordering for the positive roots
when we move from $\delta\!-\!\alpha$ to $\alpha$.
In accordance with the normal ordering (\ref{CW1}) we set

\begin{eqnarray}
\hspace{-0.7truecm}
e_{\delta}^{}\!\!\!&:=\!\!\!& [e_{\alpha}^{},e_{\delta-\alpha}^{}]_q~,
\qquad \qquad\qquad \ \ \
  e_{-\delta}^{} :=
[e_{-\delta+\alpha}^{},e_{-\alpha}^{}]_{q^{-1}}\, ,
\lb{CW3}
\\[7pt]
\hspace{-0.7truecm}
e_{n\delta+\alpha}^{}\!\!\!&:=\!\!\!& \mbox{\ls$\frac{1}{a}$}
 [e_{(n-1)\delta+\alpha}^{},e_{\delta}^{}]~,
 \qquad \ \ \;\,
 e_{-n\delta-\alpha}^{} :=  \mbox{\ls$\frac{1}{a}$}
[e_{-\delta}^{},e_{-(n-1)\delta-\alpha}^{}]\, ,
\lb{CW4}
\\[7pt]
\hspace{-0.7truecm}
e_{(n+1)\delta-\alpha}^{}\!\!\!&:=\!\!\! & \mbox{\ls$\frac{1}{a}$}
[e_{\delta}^{},e_{n\delta-\alpha}^{}]~,\;\,
\qquad \ \
 e_{-(n+1)\delta+\alpha}^{} :=  \mbox{\ls$\frac{1}{a}$}
[e_{-n\delta+\alpha}^{},e_{-\delta}^{}]\, ,
\lb{CW5}
\\[7pt]
\hspace{-0.7truecm}
e_{n\delta}'\!\!\!&:=\!\!\! & [e_{\alpha}^{},e_{n\delta-\alpha}^{}]_{q}~,
\qquad    \qquad \qquad
  e_{-n\delta}' :=
[e_{-n\delta+\alpha}^{},e_{-\alpha}^{}]_{q^{-1}}\, ,
\lb{CW6}
\edeqn
where $n=1,2,\ldots$, and $a$ is given by the formula: 
\begin{equation}
a:=[(\alpha,\alpha)]_q=\mbox{\ls$\frac{q^{(\alpha,\alpha)}-
q^{-(\alpha,\alpha)}}{q-q^{-1}}$}~.
\lb{DR14}
\end{equation}
Analogously for the inverse normal ordering (\ref{CW2}) we set
\bneqn
\hspace{-0.9truecm}
\tl{e}_{\delta}^{}\!\!\!&:=\!\!\!&
[e_{\delta-\alpha}^{},e_{\alpha}^{}]_q~,\qquad\qquad \qquad
\tl{e}_{-\delta}^{}:=[e_{-\alpha}^{},e_{-\delta+\alpha}^{}]_{q^{-1}},
\lb{CW7}
\\[7pt]
\hspace{-0.9truecm}
\tl{e}_{(n+1)\delta-\alpha}^{}\!\!\!&:=\!\!\!&\mbox{\ls$\frac{1}{a}$}
[\tl{e}_{n\delta-\alpha}^{},\tl{e}_{\delta}^{}]~,\qquad\;\,
\tl{e}_{-(n+1)\delta+\alpha}^{}:=\mbox{\ls$\frac{1}{a}$}
[\tl{e}_{-\delta}^{},\tl{e}_{-n\delta+\alpha}^{}]~,
\lb{CW8}
\\[7pt]
\hspace{-0.9truecm}
\tl{e}_{n\delta+\alpha}^{}\!\!\!&:=\!\!\!&\mbox{\ls$\frac{1}{a}$}
[\tl{e}_{\delta}^{},\tl{e}_{(n-1)\delta+\alpha}^{}]~,\qquad\;\,
\tl{e}_{-n\delta-\alpha}^{}:=\mbox{\ls$\frac{1}{a}$}
[\tl{e}_{-(n-1)\delta-\alpha}^{},\tl{e}_{-\delta}^{}]~,
\lb{CW9}
\\[7pt]
\hspace{-0.9truecm}
\tl{e}_{n\delta}'\!\!\!&:=\!\!\!&
[e_{\delta-\alpha}^{},\tl{e}_{(n-1)\delta+\alpha}^{}]_{q}~,\qquad\;\;
\tl{e}_{-n\delta}':=[e_{-\delta+\alpha}^{},
\tl{e}_{-(n-1)\delta-\alpha}^{}]_{q^{-1}},
\lb{CW10}
\edeqn
where $n=1,2,\ldots$.
Thus, we have two systems of the Cartan-Weyl generators: 'direct'
and 'inverse'. Each such system together with the Cartan generators
$k_{\alpha}^{\pm1}$, $k_{\delta-\alpha}^{\pm1}$ 
$e_{\pm\alpha}$ and $e_{\pm(\delta-a)}$ are called the q-analog of
the Cartan-Weyl basis (or simply the Cartan-Weyl basis) for
the quantum (super)algebras $U_q(A_1^{(1)})$ and $U_q(C(2)^{(2)})$.

Now we consider some properties of these bases.
First of all, the explicit construction of the Cartan-Weyl generators
(\ref{CW3})-(\ref{CW6}) (or (\ref{CW7})-(\ref{CW10})) permits easily 
to find their properties with respect to the (anti)involutions
(\ref{DR16})-(\ref{DR18}). For example, it is evident that
\begin{equation}
(e_{\pm\gamma})^*=e_{\mp\gamma}~,\qquad \forall \gamma\,\in \uDp~.
\lb{CW11}
\end{equation}
Further, it is easy to see that the 'direct' and 'inverse' Cartan-Weyl
generators ({\ref{CW3})-(\ref{CW6}) and ({\ref{CW7})-(\ref{CW10})
have very simple connection with the Dynkin involution $\tau$:
\begin{equation}
\begin{array}{rccccl}
\tau(e_{n\delta+\alpha})\!\!\!&=\!\!\!&\tl{e}_{(n+1)\delta-\alpha},\quad\;
\tau(\tl{e}_{n\delta+\alpha})\!\!\!&
=\!\!\!&e_{(n+1)\delta-\alpha}\quad &(n\in\ZZ),\quad
\\[7pt]
\tau(e_{n\delta-\alpha})\!\!\!&=\!\!\!&\tl{e}_{(n-1)\delta+\alpha},\quad\;
\tau(\tl{e}_{n\delta-\alpha})\!\!&=\!\!&e_{(n-1)\delta+\alpha}\quad &(n\in\ZZ),
\\[7pt]
\tau(e_{n\delta})\!\!\!&=\!\!\!&\tl{e}_{n\delta},\qquad\qquad\quad
\tau(\tl{e}_{n\delta})\!\!\!&=\!\!\!&e_{n\delta}\qquad\qquad &(n\neq0).
\lb{CW12}
\end{array}
\end{equation}
\setcounter{proposition}{0}
\bnpr
The root vectors (\ref{CW3})-(\ref{CW6}) satisfy the following
permutation relations:
\begin{equation}
\begin{array}{rcccl}
k_{\rm d}^{}e_{n\delta\pm\alpha}^{}k_{\rm d}^{-1}\!\!&=\!\!&
q^{n({\rm d},\delta)}e_{n\delta\pm\alpha}^{},\qquad
k_{\rm d}^{}e_{n\delta}'k_{\rm d}^{-1}\!\!&=\!\!&
q^{n({\rm d},\delta)}e_{n\delta}',
\\[7pt]
k_{\gamma}^{}e_{n\delta\pm\alpha}^{}k_{\gamma}^{}\!\!\!&=\!\!\!&
q^{\pm(\gamma,\alpha)}e_{n\delta\pm\alpha}^{},\qquad
k_{\gamma}^{}e_{n\delta}'k_{\gamma}^{-1}\!\!&=\!\!&e_{n\delta}'
\lb{CW17}
\end{array}
\end{equation}
for any $n\,\in\ZZ$ and any $\gamma\,\in\uDp$, and also
\bneqn
[e_{n\delta+\alpha},e_{-n\delta-\alpha}]&\!\!=\!\!&(-1)^{n\theta(\alpha)}
\mbox{\ls$\frac{k_{n\delta+\alpha}^{}-k_{n\delta+\alpha}^{-1}}{q-q^{-1}}$}
\quad\quad\; (n\ge0)\, ,
\lb{CW18}
\\[7pt]
[e_{n\delta-\alpha},e_{-n\delta+\alpha}]&\!\!=\!\!&(-1)^{(n-1)\theta(\alpha)}
\mbox{\ls$\frac{k_{n\delta-\alpha}^{}-k_{n\delta-\alpha}^{-1}}{q-q^{-1}}$}
\quad (n>0)\, ,
\lb{CW19}
\\[7pt]
[e_{n\delta+\alpha}^{},e_{(n+2m-\!1)\delta+\alpha}^{}]_{q}&\!\!=\!\!&
(q_{\alpha}^2\!-\!1)\!\sum_{l=1}^{m-1}q_{\alpha}^{-l}
e_{(n+l)\delta+\alpha}^{}e_{(n+2m-\!1\!-l)\delta+\alpha}^{},
\lb{CW20}
\\[7pt]
[e_{n\delta+\alpha}^{},e_{(n+2m)\delta+\alpha}^{}]_{q}&\!\!=\!\!&
(q_{\alpha}\!-\!1)q_{\alpha}^{-m+1}e_{(n+m)\delta+\alpha}^2
\nonumber \\[7pt]
&&+\,(q_{\alpha}^2\!-\!1)\!\sum\limits_{l=1}^{m-1}q_{\alpha}^{-l}
e_{(n+l)\delta+\alpha}^{}e_{(n+2m-l)\delta+\alpha}^{}\quad
\lb{CW21}
\edeqn
for any integers $n\geq0,\;m>0$;
\begin{equation}
[e_{(n+2m-1)\delta-\alpha}^{},e_{n\delta-\alpha}^{}]_{q}\!=\!
-(q_{\alpha}^2\!-\!1)\!\sum_{l=1}^{m-1}q_{\alpha}^{-l}
e_{(n+2m-\!1\!-l)\delta-\alpha}^{}e_{(n+l)\delta-\alpha}^{},\qquad
\lb{CW22}
\end{equation}
\begin{equation}
\begin{array}{rcl}
[e_{(n+2m)\delta-\alpha}^{},e_{n\delta-\alpha}^{}]_{q}&\!\!=\!\!&
-(q_{\alpha}\!-\!1)q_{\alpha}^{-m+\!1}e_{(n+m)\delta-\alpha}^2
\\[7pt]
&&-\,(q_{\alpha}^2\!-\!1)\!\sum\limits_{l=1}^{m-1}q_{\alpha}^{-l}
\!e_{(n+l)\delta-\alpha}^{}e_{(n+2m-l)\delta-\alpha}^{} \qquad
\lb{CW23}
\end{array}
\end{equation}
for any integers $n,\;m>0$;
\begin{equation}
\begin{array}{rcl}
[e_{-n\delta+\alpha}^{},e_{(n+2m-\!1)\delta+\alpha}^{}]&\!\!=\!\!&
-(\!-\!1)^{(n-1)\theta(\alpha)}(q_{\alpha}^2\!-\!1)\,
\\[7pt]
&&\!\!\!\!\!\!\!\!\!\!\!\!\!\!\!\!\!\!\!\!\!\!\!\!\!
\times\sum\limits_{l=n}^{n+m-1}\!\!q_{\alpha}^{-l}k_{n\delta-\alpha}
e_{(l-n)\delta+\alpha}^{}e_{(n+2m-\!1\!-l)\delta+\alpha}^{}
\\[7pt]
&&\!\!\!\!\!\!\!\!\!\!\!\!\!\!\!\!\!\!\!\!\!\!\!\!\!
+(q_{\alpha}^2\!-\!1)\sum\limits_{l=1}^{n-1}\!
(\!-\!1)^{l\theta(\alpha)} q_{\alpha}^{-l}
k_{\delta}^{l}e_{(-n+l)\delta+\alpha}^{}
e_{(n+2m-\!1\!-l)\delta+\alpha}^{},
\end{array}
\lb{CW24}
\end{equation}
\begin{equation}
\begin{array}{rcl}
[e_{-n\delta+\alpha}^{},e_{(n+2m)\delta+\alpha}^{}]&\!\!=\!\!&
-(-1)^{(n-1)\theta(\alpha)}(q_{\alpha}^2\!-\!1)\,
\\[7pt]
&&\!\!\!\!\!\!\!\!\!\!\!\!\!\!\!\!\!\!\!\!\!\!\!\!\!
\times\sum\limits_{l=n}^{n+m-1}\!q_{\alpha}^{-l}
k_{n\delta-\alpha}e_{(l-n)\delta+\alpha}^{}
e_{(n+2m-l)\delta+\alpha}^{}
\\[7pt]
&&\!\!\!\!\!\!\!\!\!\!\!\!\!\!\!\!\!\!\!\!\!\!\!\!\!
+(q_{\alpha}^2\!-\!1)\!\sum\limits_{l=1}^{n-1}\!
(\!-\!1)^{l\theta(\alpha)}q_{\alpha}^{-l}
k_{\delta}^{l}e_{(-n+l)\delta+\alpha}^{}e_{(n+2m-l)\delta+\alpha}^{}
\\[14pt]
&&\!\!\!\!\!\!\!\!\!\!\!\!\!\!\!\!\!\!\!\!\!\!\!\!\!
-(-1)^{(n-1)\th(\alpha)}(q_{\alpha}\!-\!1)
q_{\alpha}^{-m-n+\!1}k_{n\delta-\alpha}e_{m\delta+\alpha}^2
\end{array}
\lb{CW25}
\end{equation}
for any integers $n,\;m>0$;
\begin{equation}
\begin{array}{rcl}
[e_{(n+2m-1)\delta-\alpha}^{},e_{-n\delta-\alpha}^{}]&\!\!=\!\!&
(-1)^{(n+1)\th(\alpha)} (q_{\alpha}^2\!-\!1)
\\[7pt]
&&\!\!\!\!\!\!\!\!\!\!\!\!\!\!\!\!\!\!\!\!\!\!\!\!\!
\times\!\sum\limits_{l=n+1}^{n+m-1}\!q_{\alpha}^{-l}
e_{(n+2m-1-l)\delta-\alpha}^{}e_{(l-n)\delta-\alpha}^{}
k_{n\delta+\alpha}^{-1}
\\[7pt]
&&\!\!\!\!\!\!\!\!\!\!\!\!\!\!\!\!\!\!\!\!\!\!\!\!\!
-(q_{\alpha}^2\!-\!1)\!\sum\limits_{l=1}^{n-1}\!
(\!-\!1)^{l\theta(\alpha)}q_{\alpha}^{-l}
e_{(n+2m-1-l)\delta-\alpha}^{}e_{(-n\!+l)\delta-\alpha}^{}
k_{\delta}^{-l},
\end{array}
\lb{CW26}
\end{equation}
\begin{equation}
\begin{array}{rcl}
[e_{(n+2m)\delta-\alpha}^{},e_{-n\delta-\alpha}^{}]&\!\!=\!\!&
(-1)^{(n+1)\th(\alpha)} (q_{\alpha}^2\!-\!1)
\\[7pt]
&&\!\!\!\!\!\!\!\!\!\!\!\!\!\!\!\!\!\!\!\!\!\!\!\!\!
\times\!\sum\limits_{l=n}^{n+m-1}\!q_{\alpha}^{-l}
e_{(n+2m-l)\delta-\alpha}^{}e_{(l-n)\delta-\alpha}^{}
k_{n\delta+\alpha}^{-1}
\\[7pt]
&&\!\!\!\!\!\!\!\!\!\!\!\!\!\!\!\!\!\!\!\!\!\!\!\!\!
-(q_{\alpha}^2\!-\!1)\!\sum\limits_{l=1}^{n-1}\!
(\!-\!1)^{l\theta(\alpha)}q_{\alpha}^{-l}
e_{(n+2m-l)\delta-\alpha}^{}e_{(-n+l)\delta-\alpha}^{}k_{\delta}^{-l}
\\[14pt]
&&\!\!\!\!\!\!\!\!\!\!\!\!\!\!\!\!\!\!\!\!\!\!\!\!\!
+(-1)^{(n+1)\th(\alpha)}(q_{\alpha}\!-\!1)
q_{\alpha}^{-m-n+\!1}e_{m\delta-\alpha}^2k_{n\delta+\alpha}^{-1}
\end{array}
\lb{CW27}
\end{equation}
for any integers $n\geq0,\;m>0$;
\bneqn
[e_{n\delta+\alpha}^{},e_{m\delta-\alpha}^{}]_{q}&\!\!\!=\!\!\!&
e_{(n+m)\delta}'
\qquad\qquad\qquad\qquad\;\;(n\geq 0,\;m>0)\, ,
\lb{CW28}
\\[9pt]
[e_{n\delta+\alpha}^{},e_{-m\delta-\alpha}^{}]&\!\!\!=\!\!\!&
-(-1)^{(m+1)\th(\alpha)}
e_{(n-m)\delta}'k_{m\delta+\alpha}^{-1}\quad\,(n>m\geq 0)\, ,
\lb{CW29}
\\[9pt]
[e_{-m\delta+\alpha}^{},e_{n\delta-\alpha}^{}]&\!\!\!=\!\!\!&
-(-1)^{(m-1)\th(\alpha)}
k_{m\delta-\alpha}e_{(n-m)\delta}'\quad\,(n>m>0)\, ,
\lb{CW30}
\\[7pt]
[e_{n\delta}',e_{m\delta}']&\!\!=\!\!&[e_{-n\delta}',e_{-m\delta}']=0
\qquad\qquad\quad(n>0,\;m>0)\, ,
\lb{CW31}
\edeqn
\begin{equation}
[e_{n\delta+\alpha}^{},e_{m\delta}']=
q_{\alpha}^{-m+1}ae_{(n+m)\delta+\alpha}^{}+
(q_{\alpha}^2\!-\!1)\!\sum\limits_{l=1}^{m-1}q_{\alpha}^{-l}
e_{(n+l)\delta+\alpha}^{}e_{(m-l)\delta}'
\lb{CW32}
\end{equation}
for any integers $n\geq0,\;m>0$;
\begin{equation}
[e_{m\delta}',e_{n\delta-\alpha}^{}]=
q_{\alpha}^{-m+\!1}ae_{(n+m)\delta-\alpha}^{}+
(q_{\alpha}^2\!-\!1)\!\sum\limits_{l=1}^{m-1}q_{\alpha}^{-l}
\!e_{(m-l)\delta}'e_{(n+l)\delta-\alpha}^{}
\lb{CW33}
\end{equation}
for any integers $n,\;m>0$;
\begin{equation}
\begin{array}{rcl}
[e_{-n\delta+\alpha}^{},e_{m\delta}']&\!\!=\!\!&-
(-1)^{(n-1)\th(\alpha)}q_{\alpha}^{-m+1}a
k_{n\delta-\alpha}^{}e_{(m-n)\delta+\alpha}^{}
\\[7pt]
&&-(-1)^{(n-1)\th(\alpha)}(q_{\alpha}^2\!-\!1)
k_{n\delta-\alpha}^{}\sum\limits_{l=n}^{m-1}\!q_{\alpha}^{-l}
e_{(l-n)\delta+\alpha}^{}e_{(m-l)\delta}'
\\[7pt]
&&+(q_{\alpha}^2\!-\!1)\!\sum\limits_{l=1}^{n-1}\!
(-1)^{l\th(\alpha)}q_{\alpha}^{-l}k_{\delta}^{l}
e_{(-n+l)\delta+\alpha}^{}e_{(m-l)\delta}'
\end{array}
\lb{CW34}
\end{equation}
for any integers $m\ge n>0$;
\begin{equation}
\begin{array}{rcl}
[e_{-n\delta+\alpha}^{},e_{m\delta}']&\!\!=\!\!&
(-1)^{m\th(\alpha)}q_{\alpha}^{-m+1}a
k_{\delta}^{m}e_{(-n+m)\delta+\alpha}^{}
\\[7pt]
&&+(q_{\alpha}^2\!-\!1)\!\sum\limits_{l=1}^{m-1}\!
(-1)^{l\th(\alpha)}q_{\alpha}^{-l}k_{\delta}^{l}
e_{(-n+l)\delta+\alpha}^{}e_{(m-l)\delta}'
\end{array}
\lb{CW35}
\end{equation}
for any integers $n>m>0$;
\begin{equation}
\begin{array}{rcl}
[e_{m\delta}'e_{-n\delta-\alpha}^{}]&\!\!=\!\!&
-(-1)^{(n+1)\th(\alpha)}q_{\alpha}^{-m+1}a
e_{(m-n)\delta-\alpha}^{}k_{n\delta+\alpha}^{-1}
\\[7pt]
&&-(-1)^{(n+1)\th(\alpha)}(q_{\alpha}^2\!-\!1)\!
\sum\limits_{l=n+1}^{m-1}\!q_{\alpha}^{-l}
e_{(m-l)\delta}'e_{(l-n)\delta-\alpha}^{}k_{n\delta+\alpha}^{-1}
\\[7pt]
&&+(q_{\alpha}^2\!-\!1)\!\sum\limits_{l=1}^{n}\!
(-1)^{l\th(\alpha)}q_{\alpha}^{-l}
e_{(m-l)\delta}'e_{(-n+l)\delta-\alpha}^{}k_{\delta}^{-l}
\end{array}
\lb{CW36}
\end{equation}
for any integers $m>n\ge0$;
\begin{equation}
\begin{array}{rcl}
[e_{m\delta}'e_{-n\delta-\alpha}^{}]&\!\!=\!\!&
(-1)^{m\th(\alpha)}q_{\alpha}^{-m+1}a
e_{(-n+m)\delta-\alpha}^{}k_{\delta}^{-m}
\\[7pt]
&&+(q_{\alpha}^2\!-\!1)\!\sum\limits_{l=1}^{m-1}\!
(-1)^{l\th(\alpha)}q_{\alpha}^{-l}
e_{(m-l)\delta}'e_{(-n+l)\delta-\alpha}^{}k_{\delta}^{-l}
\end{array}
\lb{CW37}
\end{equation}
for any integers $n\ge m>0$.
\lb{PCW1}
\edpr
Here in the relations (\ref{CW20})-(\ref{CW37}) and in what follows
we denote $q_{\alpha}\!:=\!(-1)^{\th(\alpha)}q^{(\alpha,\alpha)}$.
The imaginary root vectors $e_{n\delta}'$ do not satisfy the relations 
of the type (\ref{CW18}) and therefore we introduce  new imaginary
roots vectors $e_{\pm n\delta}$ by the following (Schur) relations:
\begin{equation}
{e'}_{n\delta}=\sum_{p_{1}+2p_{2}+\ldots+np_{n}=n}\!\!\!\!\!
\mbox{\ls$\frac{\;\Bigl((-1)^{\th(\alpha)}
(q-q^{-1})\Bigr)^{\sum p_{i}-1}}{p_{1}!\cdots
p_{n}!}$}\;\;e_{\delta}^{\;p_1}\cdots e_{n\delta}^{\;p_n}.
\lb{CW38}
\end{equation}
In terms of generating functions
\bneqn
{\cal E}'(u)&\!\!:=\!\!&(q-q^{-1})\sum_{n\geq 1}^{}{e'}_{n\delta}u^{-n}~,
\lb{CW39}
\\[7pt]
{\cal E}(u)&\!\!=\!\!&(q-q^{-1})\sum_{n\geq 1}^{}e_{n\delta}u^{-n}
\lb{CW40}
\edeqn
the relation (\ref{CW38}) may be rewritten in the form
\begin{equation}
{\cal E}'(u)=-1+\exp{\cal E}(u)
\lb{CW41}
\end{equation}
or
\begin{equation}
{\cal E}(u)=\ln(1+{\cal E}'(u))~.
\lb{CW42}
\end{equation}
This provides a formula inverse to (\ref{CW38})
\begin{equation}
e_{n\delta}=\!\!\!\sum_{p_{1}+2p_{2}+\ldots+np_{n}=n}\!\!\!\!\!\!\!\!\!\!
\!\!\!\mbox{\ls$\frac{\Bigl((-1)^{\th(\alpha)}
(q^{-1}-q)\Bigr)^{\sum p_{i}-1}(\sum_{i=1}^{n}p_{i}-1)!}
{p_{1}!\cdots p_{n}!}$}(e_{\delta}')^{p_1}\cdots 
(e_{n\delta}')^{p_n}\!\!.
\lb{CW43}
\end{equation}
The new root vectors corresponding to negative roots are obtained 
by the Cartan conjugation $(^*)$:
\begin{equation}
 e_{-n\delta}=(e_{n\delta})^{*}~.
\lb{CW44}
\end{equation}
\bnpr
The new root vectors $e_{\pm n\delta}$ satisfy the following
commutation relations:
\bneqn
[e_{n\delta+\alpha}^{},e_{m\delta}^{}]&\!\!\!=\!\!\!&
(\!-1\!)^{(m-1)\th(\alpha)}a(\!m\!)e_{(n+m)\delta+\alpha}^{}
\qquad\;(n\ge0,\,m>0),
\lb{CW45}
\\[7pt]
[e_{m\delta}^{},e_{n\delta-\alpha}^{}]&\!\!\!=\!\!\!&
(\!-1\!)^{(m-1)\th(\alpha)}a(\!m\!)e_{(n+m)\delta-\alpha}^{}
\qquad\quad\quad\;(n,\,m>0),
\lb{CW46}
\\[7pt]
[e_{-n\delta+\alpha}^{},e_{m\delta}^{}]&\!\!\!=\!\!\!&
-(\!-1\!)^{(n+m)\th(\alpha)}a(\!m\!)
k_{n\delta-\alpha}^{}e_{(m-n)\delta+\alpha}^{}\;\;
(m\ge n>0),\phantom{pls} 
\lb{CW47}
\\[7pt]
[e_{-n\delta+\alpha}^{},e_{m\delta}^{}]&\!\!\!=\!\!\!&
(\!-1\!)^{\th(\alpha)}\!a(\!m\!)k_{\delta}^{m}
e_{(-n+m)\delta+\alpha}^{}\qquad\qquad\,(n>m>0),
\lb{CW48}
\\[7pt]
[e_{m\delta}^{},e_{-n\delta-\alpha}^{}]&\!\!\!=\!\!\!&
-(\!-1\!)^{(n+m)\th(\alpha)}\!a(\!m\!)e_{(m-n)\delta-\alpha}^{}
k_{n\delta+\alpha}^{-1}\;\;\,(m>n\ge0),
\lb{CW49}
\\[7pt]
[e_{n\delta-\alpha}^{},e_{-m\delta}^{}]&\!\!\!=\!\!\!&
(\!-1\!)^{\th(\alpha)}\!a(\!m\!)e_{(-n+m)\delta-\alpha}^{}
k_{\delta}^{-m}\qquad\quad\;\;\,(n\ge m>0),
\lb{CW50}
\\[7pt]
[e_{n\delta},e_{-m\delta}]&\!\!\!=\!\!\!&\delta_{nm}\,a(m)\,
\mbox{\ls$\frac{k_{\delta}^{m}-k_{\delta}^{-m}}{q-q^{-1}}$}
\qquad\qquad\qquad\qquad(n,\,m>0),
\lb{CW51}
\edeqn
where
\begin{equation}
a(m):=\mbox{\ls$\frac{q^{m(\alpha,\alpha)}-q^{-m(\alpha,\alpha)}}
{m(q-q^{-1})}$}~.
\lb{CW52}
\end{equation}
\lb{PCW2}
\edpr
All the relations of Propositions (\ref{PCW1}), (\ref{PCW2}) together
with the ones obtained from them by the conjugation describe complete
list of the permutation relations of the Cartan-Weyl bases
corresponding to the 'direct' normal ordering (\ref{CW1}).
Applying to these relations the Dynkin involution $\tau$, it is easy
to obtain these results for the 'inverse' normal ordering
(\ref{CW2}).

\setcounter{equation}{0}
\section{Universal $R$-matrix for $\!U_{\!q}(\!A_1^{(\!1\!)})$ and
$\!U_{\!q}(\!C(2)^{(\!2\!)})$}

Any quantum (super)algebra $U_{q}(g)$ is a non-cocommutative Hopf
(super)\-al\-gebra which has the intertwining operator called the
universal $R$-matrix.
By definition \cite{D}, the universal $R$-matrix for the Hopf
(super)algebra $U_q(g)$ is an invertible element $R$ of some
extension $U_q(g)\ot U_q(g)$, satisfying the equations
\bneqn
\tl{\Delta}_q(a)&\!\!=\!\!&R\Delta_q(a)R^{-1} \qquad\quad\;\;
\forall\,\,a \in U_q(g)~,
\lb{UR1}
\\[7pt]
(\Delta_q\ot\id)R&\!\!=\!\!&R^{13}R^{23}~,\qquad
(\id\ot\Delta_q)R=R^{13}R^{12},
\lb{UR2}
\edeqn
where  $\tl{\Delta}_q$ is the opposite comultiplication:
\sloppy $\tl{\Delta}_q=\sg\Delta_q$, 
$\sg(a\ot b)=(-1)^{\deg a\deg b} b\ot a$ for
all homogeneous elements $a,\,b \in U_q(g)$. In the relation
(\ref{UR2}) we use the standard notations
$R^{12}=\sum a_{i}\ot b_{i}\ot\id$,
$R^{13}=\sum a_{i}\ot\id\ot b_{i}$,
$R^{23}=\sum \id\ot a_{i}\ot b_{i}$ if $R$ has the form
$R=\sum a_{i}\ot b_{i}$.
We employ the following standard notation for the q-exponential:
\begin{equation}
\exp_{q}(x):=1+x+\mbox{\ls$\frac{x^{2}}{(2)_{q}!}$}+\ld+
\mbox{\ls$\frac{x^{n}}{(n)_{q}!}$}+\ld=\sum_{n\geq0}
\mbox{\ls$\frac{x^{n}}{(n)_{q}!}$}~,
\lb{UR3}
\end{equation}
where
\begin{equation}
(n)_{q}:=\mbox{\ls$\frac{q^{n}-1}{q-1}$}~.
\lb{UR4}
\end{equation}
An explicit expression of the universal $R$-matrix $R$ for our case
$U_q(g)$ can be presented as follows:
\begin{equation}
R=R_{+}R_{0}R_{-}K~.
\lb{UR5}
\end{equation}
Here the factors $K$ and $R_{\pm}$ have the following form
\bneqn
K&\!\!=\!\!&q^{\frac{1}{(\alpha,\alpha)}h_{\alpha}\ot h_{\alpha}+
h_{\delta}\ot h_{\rm d}+h_{\rm d}\ot h_{\delta}}~,
\lb{UR6}
\\[9pt]
R_{+}&\!\!=\!\!&\prod_{n\ge0}^{\ra}\cR_{n\delta+\alpha}, \qquad
R_{-}=\prod_{n\ge1}^{\la}\cR_{n\delta-\alpha}.
\lb{UR7}
\edeqn
were the elements $\cR_{\gamma}$ are given by the formula
\begin{equation}
\cR_{\gamma}=\exp_{q^{-1}_{\gamma}}\Bigl(A(\gamma)
(q-q^{-1})(e_{\gamma}\ot e_{-\gamma})\Bigr)~,
\lb{UR8}
\end{equation}
where
\bneqn
q_{\gamma}&\!\!=\!\!&(-1)^{\vth(\gamma)}q^{(\gamma,\gamma)}~,
\lb{UR9}
\\[9pt]
A(\gamma)&\!\!=\!\!&\left\{\begin{array}{lll}
&\!\!\!\!\!(-1)^{n\th(\alpha)}\quad &{\rm if}\;\gamma=n\delta+\alpha~,\\
&\!\!\!\!\!(-1)^{(n-1)\th(\alpha)} \quad &{\rm if}\;\gamma=n\delta-\alpha~.
\end{array}\right.
\label{UR10}
\edeqn
Finally, the factor $R_{0}$ is defined as follows
\begin{equation}
R_{0}=\exp\Bigl((q-q^{-1})\sum_{n>0}^{}d(n)
e_{n\delta}\ot e_{-n\delta}\Bigr)~,
\label{UR11}
\end{equation}
where $d(n)$ is the inverse to $a(n)$, i.e.
\begin{equation}
d(n)=\mbox{\ls$\frac{n(q-q^{-1})}
{q^{n(\alpha,\alpha)}-q^{-n(\alpha,\alpha)}}$}~.
\lb{UR12}
\end{equation}
\section*{Acknowledgments}
The third author is grateful to Institute of Theoretical Physics,
University of Wroclaw and the Organizing Committee of XII Max Born
Symposium for the support of his visit on the Symposium.
S.M. Khoroshkin, V.N. Tolstoy  would like to thank for the
support by the Russian Foundation for Fundamental Research (grant
No.98-01-00303) as well as
 by the program of French-Russian scientific cooperation
(CNRS grant PICS-608 and  grant RFBR-98-01-22033)  
and J. Lukierski for the support by KBN grant 2PO3B13012.


\begin{thebibliography}{99}
\bibitem{AST}
Asherova, R.M., Smirnov, Yu.F., and Tolstoy, V.N.
A description of certain class of projection operators for
complex semisimple Lie algebras. (Russian) 
{\it Matem. Zametki} {\bf 26} (1979), no. 3, 15--25.

\bibitem{D}
Drinfeld, V.G.
Quantum groups.
{\it Proc. ICM-86 (Berkely USA) vol.1}, 798--820. Amer. Math. Soc.
(1987).

\bibitem{K1}
Kac, V.G.
Lie superalgebras.
{\it Adv. Math.} {\bf  26} (1977), 8--96.

\bibitem{K2}
Kac, V.G.
Infinite dimensional Lie  algebras.
{\it Cambridge University Press, 1990}.

\bibitem{KT1}
Khoroshkin, S.M., and Tolstoy, V.N.
Universal R-matrix for quan\-tized (super)\-algebras.
{\it Commun. Math. Phys.} {\bf 141} (1991), no. 3, 599--617.

\bibitem{KT2}
Khoroshkin, S.M., and Tolstoy, V.N.
The Cartan-Weyl basis and the universal $R$-matrix for quantum
Kac-Moody algebras and superalgebras.
{\it Quantum Symmetries (Clausthal 1991)}. 336--351,
{\it World Sci. Publishing, River Edge, NJ}, 1993.

\bibitem{KT3}
Khoroshkin, S.M., and Tolstoy, V.N.
The uniqueness theorem for the universal R-matrix.
{\it Lett. Math. Phys.} {\bf  24} (1992), no. 3, 231--244.

\bibitem{KT4}
Khoroshkin S.M., and Tolstoy V.N.
Extremal projector and universal $R$-matrix for quantum
contragredient Lie (super)algebras.
{\it Quantum groups and related topics (Wroclaw, 1991)}, 23--32,
Math. Phys. Stud., 13, {\it Kluwer Acad. Publ., Dordrecht}, 1992.

\bibitem{nieznany}
Lukierski, J., and Tolstoy, V.N. 
Cartan-Weyl basis for quantum affine superalgebra 
$U_q(\widehat{\rm osp}\,(1|2))$.
{\it Czech. Journ. Phys.} {\bf 47} (1997), no. 12, 1231--1239.

\bibitem{L}
Lusztig, G.
Canonical bases arising from quantized enveloping algebras.
{\it J. Amer. Math. Soc.} {\bf 3} (1990), 447--498.

\bibitem{T1}
Tolstoy, V.N.
Extremal projectors for quantized Kac-Moody superalgebras and some
of their applications. {\it Quantum Groups (Clausthal, 1989)},
118--125, Lectures Notes in Phys. 370, {\it Springer, Berlin}, 1990.

\bibitem{T2}
Tolstoy, V.N.
Extremal projectors for contragredient Lie algebras and
superalgebras of finite growth. (Russian) {\it Uspekhi Math. Nauk}
{\bf 44} (1989), no.1(265), 211--212; {\it translation in Russian 
Math. Serveys} {\bf 44} (1989), no. 1, 257--258. 

\bibitem{TK}
Tolstoy, V.N., and Khoroshkin, S.M.
The Universal R-matrix for quantum nontwisted affine Lie algebras.
(Russian) {\it Funktsional. Anal. i Prilozhen.}  {\bf  26} (1992), 
no. 1, 85--88; {\it translation in Functional Anal. Appl.} 
{\bf 26} (1992), no. 1, 69--71.

\bibitem{VdL}
Van der Leur, J.W. Contragredient Lie superalgebras of finite grouth.
{\it Utrecht thesis}, 1985. 
\end{thebibliography}
\end{document}